\documentclass[a4paper,11pt]{article}

\usepackage{array}
\usepackage{theorem}
\usepackage{amsmath,amscd,amssymb} 
\usepackage{tikz}
\usepackage{latexsym}
\usepackage{enumerate}

\theorembodyfont{\sl}

\newtheorem{lemma}{Lemma}[section]

\newtheorem{proposition}[lemma]{Proposition}
\newtheorem{theorem}[lemma]{Theorem}
\newtheorem{corollary}[lemma]{Corollary}

\newtheorem{definition}[lemma]{Definition}

\newtheorem{introthm}{Theorem}

\renewcommand{\AA}{\mathbb A}

\newcommand{\CC}{\mathbb C}

\newcommand{\HH}{\mathbb H}

\newcommand{\QQ}{\mathbb Q}
\newcommand{\RR}{\mathbb R}

\newcommand{\TT}{\mathbb T}

\newcommand{\WW}{\mathbb W}
\newcommand{\XX}{\mathbb X}

\newcommand{\ZZ}{\mathbb Z}

\newcommand{\cA}{\mathcal A}

\newcommand{\cD}{\mathcal D}

\newcommand{\cM}{\mathcal M}

\newcommand{\cO}{\mathcal O}

\renewcommand{\Bar}{\overline}
\newcommand{\half}{{\frac{1}{2}}}

\newcommand{\imic}{\cong}

\newcommand{\splitextn}{\ltimes}

\newcommand{\tensor}{\otimes}
\renewcommand{\Tilde}{\widetilde}

\newcommand{\GL}{\mathop{\mathrm {GL}}\nolimits}
\newcommand{\Sp}{\mathop{\mathrm {Sp}}\nolimits}
\newcommand{\Ut}{\mathop{\null\mathrm {U}}\nolimits}
\newcommand{\Gr}{\mathop{\mathrm {Gr}}\nolimits}
\newcommand{\Pic}{\mathop{\mathrm {Pic}}\nolimits}
\newcommand{\Sym}{\mathop{\mathrm {Sym}}\nolimits}
\newcommand{\Imag}{\mathop{\mathrm {Im}}\nolimits}

\newcommand{\rank}{\mathop{\mathrm {rank}}\nolimits}
\newcommand{\rad}{\mathop{\mathrm {rad}}\nolimits}
\newcommand{\reg}{\mathop{\mathrm {reg}}\nolimits}
\newcommand{\tr}{\mathop{\mathrm {tr}}\nolimits}
\newcommand{\Eins}{{\mathbf 1}}

\newcommand{\latt}[1]{{\langle{#1}\rangle}}
\newcommand{\Null}{\mathop{\mathrm {null}}\nolimits}
\newcommand{\sym}{\mathop{\mathrm {sym}}\nolimits}
\newcommand{\Temb}{\mathop{\mathrm {Temb}}\nolimits}
\newcommand{\RT}{\mathop{\mathrm {RT}}\nolimits}
\newcommand{\qed}{\hfill$\blacksquare$}

\begin{document}
\title{Slopes of Siegel cusp forms and geometry of compactified Kuga varieties}
 \author{Flora Poon, Riccardo Salvati Manni and Gregory Sankaran}
%
%
%
%
%\author{Flora Poon}
%\email{wkpoon@ncts.ntu.edu.tw\,\,\, ORCID:0009-0005-0343-6920}
%\address{Mathematics Division, National
%  Center for Theoretical Sciences, \newline
%  No.~1, Sec.~4, Roosevelt~Rd., Taipei City, Taiwan.}
%
%\author{Riccardo Salvati Manni}
%\email{salvati@mat.uniroma1.it\,\,\,  ORCID:0000-0001-7467-2625}
%\address{Dipartimento di Matematica ``Guido Castelnuovo'', Universit\`a di Roma ``La Sapienza'',\newline P.le~Aldo~Moro~5, 00185~Roma, Italy}
%
%\author{Gregory Sankaran}
%\email{G.K.Sankaran@bath.ac.uk,\,\,  ORCID: 0000-0002-5846-6490}
%\address{Department of Mathematical Sciences, University
%  of Bath,\newline Bath~BA2~7AY, UK}
%
%\classification{14K10 (primary), 11F46 (secondary).  \,\,\,\, \quad\quad
% DOI:10.48550/arxiv.2109.06142}
%\keywords{Kodaira dimension, Kuga variety, Siegel modular forms, Toroidal compactification.}

%\begin{abstract}
%We study the Kodaira dimension of the compactified $n$-fold Kuga
%variety over the moduli space of principally polarised abelian
%$g$-folds. We construct a suitable compactification, which we call a
%Namikawa compactification, and show that in most cases it has
%canonical singularities. We then use results about the slope of Siegel
%modular forms to determine the Kodaira dimension for all $g>1$ and
%$n>0$.
%\end{abstract}

\maketitle

\section{Introduction}\label{sect:intro}

The universal family $X(\Gamma)$ over a moduli space 
$\cA(\Gamma)=\Gamma\backslash\HH_g$ of abelian varieties associated with a finite
index subgroup of $\Sp(2g,\ZZ)$ is known as the Kuga variety. Such
families were first studied systematically by M.~Kuga, whose 1964
Chicago lecture notes on the subject~\cite{kuga} have been recently
published. The construction is given in~\cite{kuga} for
$\Gamma<\Sp(2g,\ZZ)$ a torsion-free subgroup of finite index, but the
restriction to torsion-free can be removed. However, if $-1\in\Gamma$
then the fibre of the family is no longer the abelian variety but
instead the corresponding Kummer variety. One could also allow
$\Gamma$ to be a subgroup of $\Sp(2g,\QQ)$ commensurable with
$\Sp(2g,\ZZ)$, for example taking $\cA(\Gamma)$ to be the moduli space
of abelian varieties with some non-principal polarisation, but we
shall not pursue this here.

A natural generalisation is to consider the $n$-fold Kuga variety
$X^n(\Gamma)$, whose general fibre is the $n$-fold direct product $A^n$ of
the corresponding abelian variety $A$ or, if $-1\in\Gamma$, the
Kummer variety $A^n/\pm 1$.

Alternatively one may consider the universal family
$\underline{X}^n(\Gamma)$ over the stack
$\underline{\cA}(\Gamma):=[\Gamma\backslash\HH_g]$. In this case the
fibre is an abelian variety in all cases, but if $-1\in \Gamma$ then
the base has non-trivial stabilisers generically. This is the object
that is studied in a particular case in~\cite{fv}.

We shall be concerned with compactifications of the $n$-fold Kuga
variety $X^n_g=X^n(\Sp(2g,\ZZ))$ associated with the coarse moduli
space $\cA_g$ of principally polarised abelian $g$-folds over
$\CC$. Thus~$X_g^1=X(\Sp(2g,\ZZ))$, and by convention $X_g^0=\cA_g$.

The starting point for our work is Ma's study~\cite{ma} of
$X_g^n(\Gamma)$. We construct a special compactification, which we
call a Namikawa compactification, of $X^n_g$ and this, together with
recent and less recent results about the slope of $\cA_g$, allows
us to determine the Kodaira dimension~$\kappa(X_g^n)$ of $X_g^n$ whenever
$g\ge 2$ and $n\ge 1$.
 
\begin{introthm}\label{thm:kdim}
Suppose that $g\ge 2$ and $n\ge 1$. Then the Kodaira dimension
$\kappa(X^n_g)$ of
$X^n_g$ satisfies:
\begin{itemize}
\item
  $\kappa(X^n_g)=\half g(g+1)=\dim\cA_g$ if $g+n\ge 7$, except for
  $(g,n)= (4,3)$, $(3,5)$, $(3,4)$, $(2,7)$, $(2,6)$ and $(2,5)$;
\item
  $\kappa(X_2^7)=\kappa(X^5_3)=\kappa(X^3_4)=0$; and
\item
 $\kappa(X^n_g)=-\infty$ otherwise, i.e.\ if $g+n\leq 6$ or
  $(g,n)=(2,6)$, $(2,5)$ or $(3,4)$.
\end{itemize}
\end{introthm}

For completeness, we mention the case $g=1$. We can identify $X_1^n$
with $\cM_{1,n+1}$, the coarse moduli space of genus~$1$ curves with
$n+1$ (ordered) marked points: the case $n=0$ is the familiar
identification of $\cM_{1,1}$ with $\AA^1$ via the $j$-invariant. The
results of Belorousski~\cite{belo} then show that $X_1^n$ is rational
when $n\leq 9$, and from~\cite{bf} we have $\kappa(X_1^{10})=0$ and
$\kappa( X_1^{n})=1$ for $n\geq 11$.

In fact $X_2^n$ for $n\le 6$ and $X^n_3$ for $n\le 4$ are all
unirational. If $g\le 3$ then, since $\cA_g$ is birationally
equivalent to $\cM_g$, and the $g$th symmetric power of a curve is
birational to its Jacobian, there is a dominating rational map
$\cM_{g,gn}\dasharrow X_g^n$. On the other hand, it is known that
$\cM_{2,n}$ and $\cM_{3,n}$ are unirational for $n<13$ and $n<15$
respectively: see \cite[Table~1]{tant}.

Moreover $X_g^1$ is unirational for $g\le 5$: see~\cite{verra} for
$g=3$ and $g=4$ and~\cite{fv} for $g=5$. Note that $X_g^1$ is the
boundary of the Mumford partial compactification $\cA'_{g+1}$. One can
therefore compactify $X_g^1$ by taking its closure in any toroidal
compactification of $\cA_{g+1}$, since these all contain $\cA'_{g+1}$
as a dense open set, but it is not straightforward to control the
singularities that then arise.
\medskip

For the cases $g+n\le 6$ our proof is based on the fact that the slope
of $\cA_g$ is greater than $7$: this was proved in \cite{sm} for $g=3$ and
$g=4$, and in~\cite{fgsv} for $g=5$.

For the cases $g+n\geq 7$, our main technical result concerns the
existence of a sufficiently good compactification of $X_g^n$. We will
say that a compactification of $X_g^n$ is a Namikawa compactification
if it dominates a toroidal compactification of $\cA_g$ and boundary
divisors are mapped to boundary divisors: see
Definition~\ref{def:namikawa} for a full explanation and precise
details. We prove (see Theorem~\ref{thm:goodcompact} for the precise
statement):

\begin{introthm}\label{thm:canonicalcompact}
Suppose $g\ge 2$ and $n\ge 1$. There exists a Namikawa
compactification $\Bar{X_g^n}$ with canonical singularities as long as
$g+n\ge 6$.
\end{introthm}

There is some overlap between our results and those of~\cite{ma}.  The
singularities of $X_g^n$ are studied in~\cite[Section~10]{ma}, but the
singularities at the boundary $\Bar{X_g^n}$ are not considered
there. That is sufficient for computing the geometric genus, but not
other plurigenera: on the other hand, Ma's approach gives precise
information about the geometric genus and therefore some information
about the Kodaira dimension.

Our original motivation came from the case $g=6$ and $n=1$, and the
Kodaira dimension~$\kappa(X_6^1)$.  As we have seen, $X_g^1$ is
unirational for $g<6$.  On the other hand, for $g\geq 7$ we have~$\kappa(X_g^1)=\kappa(\cA_g)= \half g(g+1)$, since Iitaka's conjecture
holds, cf.~\cite{kaw}, and $\cA_g$ is of general type.

Similar questions arise in relation to the universal Jacobian
varieties over $\cM_g$. A suitable compactification was constructed by
Caporaso~\cite{cap} and the question of Kodaira dimension is well
studied: see for example \cite{bfv, cmkv, fvjac}, as well as \cite{bf}
for the case $g=1$.

\section{Namikawa compactifications}\label{sect:namikawa}

Suppose throughout this section that $g\ge 2$. Let $\cA_g\subset
\Bar\cA_g$ be the inclusion of the coarse moduli space
$\cA_g=\Sp(2g,\ZZ)\backslash\HH_g$ of principally polarised abelian varieties
of dimension $g$ into a toroidal compactification $\Bar\cA_g$. Denote
by $f\colon X_g^n=(\ZZ^{2gn}\splitextn\Sp(2g,\ZZ))\backslash(\CC^{gn}\times \HH_g)$
the $n$-fold Kuga family, as in~\cite{ma}.

\begin{definition}\label{def:namikawa}
A \emph{Namikawa compactification} of $X_g^n$ is an irreducible normal
projective variety~$\Bar{X_g^n}$ containing $X_g^n$ as an open subset,
together with a projective toroidal compactification~$\Bar\cA_g$ of~$\cA_g$ for which the following conditions hold.
\begin{enumerate}
  \item $f\colon X_g^n\to \cA_g$ extends to a projective morphism
    $f\colon\Bar{X_g^n}\to \Bar\cA_g$;
  \item every irreducible component of
    $\Delta_X:=\Bar{X_g^n}\smallsetminus X_g^n$ dominates an
    irreducible component of~${\Delta_\cA:=\Bar\cA_g\smallsetminus\cA_g}$; 
\end{enumerate}
\end{definition}

Compactifications satisfying these conditions were first given by
Na\-mi\-ka\-wa~\cite{nam,nambook}, for the case $n=1$. There is little
difficulty in extending Namikawa's construction to arbitrary~$n$, and it is
essentially done in~\cite[Ch.~VI.1]{fc}. The conditions in
Definition~\ref{def:namikawa} are the same as in
\cite[Theorem~1.2]{ma}, but we also require $\Bar{X_g^n}$ to be normal
(rather than just smooth in codimension~$1$) and projective. In fact
neither of these conditions presents any difficulty.

\begin{theorem}\label{thm:goodcompact}
Suppose $g\ge 2$ and $n\ge 1$. Then if $g+n\ge 6$, there exists a
Namikawa compactification $\Bar{X_g^n}\supset X_g^n$ such that
$\Bar{X_g^n}$ has canonical singularities. In particular $\Bar{X_g^n}$
is $\QQ$-Gorenstein.
\end{theorem}

The proof will occupy the rest of this section. The first step is the
following lemma. Although similar statements are already known (see
for example~\cite{as-b}) we give a proof here as we have not found one
elsewhere.

Suppose that a finite group $G$ acts effectively on a variety $X$. A
nontrivial element of $G$ is called a \emph{quasireflection} if it
preserves a divisor on $X$.

\begin{lemma}\label{lem:covercanonical}
Suppose that $G$ is a finite group acting effectively and without
quasireflections on a variety $X$ that has canonical
singularities. Let $f\colon \hat X\to X$ be a $G$-equivariant
resolution of singularities, and suppose that $\hat X/G$ has canonical
singularities. Then $X/G$ also has canonical singularities.
\end{lemma}

\noindent{\it Proof.\/} Since $X$ has canonical singularities it is in
particular $\QQ$-Gorenstein (we do not require, nor expect, $X$ or
$X/G$ to be $\QQ$-factorial), and therefore $X/G$ is also
$\QQ$-Gorenstein. Suppose that~$rK_{X/G}$ is Cartier and $\sigma\in
\cO(rK_{X/G})$, so $\sigma$ is a pluricanonical form on
$(X/G)_{\reg}$. Therefore~$\sigma$ lifts to a $G$-invariant form
$g^*\sigma$ on an open $G$-invariant subset $X_0\subset X$, where
$g\colon X\to X/G$ is the quotient map.

The complement $X\smallsetminus X_0$ consists of the fixed loci of
elements of $G$, together with the singular locus of $X$; but the
fixed loci have codimension at least~$2$ by assumption, so $g^*\sigma$
extends $G$-invariantly to $X_{\reg}$. Therefore it lifts on $\hat X$
to a form defined away from the exceptional locus of $f$, but because
$X$ has canonical singularities, this extends to a $G$-invariant form~$\hat\sigma = f^*g^*\sigma$ on $\hat X$ without poles, which in turns
descends to the smooth part $(\hat X/G)_{\reg}$ of $\hat X/G$, because~${\hat g \colon \hat X\to \hat X/G}$ is \'etale in codimension~$1$.

This form $\hat\sigma$ agrees with $\bar f^*\sigma$ on the dense open
set $\bar f^{-1}((X/G)_{\reg})$, where $\bar f\colon \hat X/G \to X/G$
is the map induced by $f$: therefore $\bar f^*\sigma$ extends to
$(\hat X/G)_{\reg}$.

Now consider a resolution of singularities $h\colon Y \to X/G$ and
form the pullback $\hat h\colon \hat Y\to \hat X/G$ as in the diagram.
\begin{center}
\begin{tikzpicture}
  \node (Xh) at (-1,2) {$\hat X$};
  \node (X) at (3,2) {$X$};
  \node (XGh) at (0,0) {$\hat X/G$};
  \node (XG) at (4,0) {$X/G$};
  \node (Yh) at (-1,-2) {$\hat Y$};
  \node (Y) at (3,-2) {$Y$};
  \node [above] at (1,2) {$f$};
  \node [left] at (-0.5,1) {$\hat g$};
  \node [left] at (3.5,1) {$g$};
  \node [above] at (2,0) {$\bar f$};
  \node [left] at (-0.5,-1) {$\hat h$};
  \node [left] at (3.5,-1) {$h$};
  \node [above] at (1,-2) {$\hat f$};
  \draw [->, >=latex](Xh.east)--(X.west);
  \draw [->, >=latex](Xh.south east)--(XGh.north);
  \draw [->, >=latex](XGh.east)--(XG.west);
  \draw [->, >=latex](X.south east)--(XG.north);
  \draw [->, >=latex](Yh.north east)--(XGh.south);
  \draw [->, >=latex](Yh.east)--(Y.west);
  \draw [->, >=latex](Y.north east)--(XG.south);
\end{tikzpicture}
\end{center}

This resolves the singularities of $\hat X/G$, so $\hat h^* \bar
f^*\sigma$ extends without poles to the whole of $\hat Y$. But now if
$h^*\sigma$ has poles along a divisor $E\subset Y$, then $\hat f^* h^*
\sigma=\hat h^* \bar f^*\sigma$ has poles along $\hat f^{-1}E$, which
is impossible. Therefore $h^*\sigma$ is holomorphic, and hence $X/G$
has canonical singularities. \qed\smallskip
\medskip

Next, we construct Namikawa compactifications $\Bar{X_g^n}$. To do
this we follow~\cite{nam} closely, and as far as possible we have
chosen our notation to be compatible with that paper. The notation
in~\cite{fc} is different. Note that our compactifications are
singular: both \cite[Proposition~5.4(i)]{nam} and
\cite[Theorem~VI.1.1(i)]{fc} mention smooth compactifications, but
after a base change in the first case and as stacks in the second (see
the remarks at the top of~\cite[Page~195]{fc}).

We write $\Gamma=\Sp(2g,\ZZ)$. For simplicity, and because it is
enough for our purposes, we consider only this case, but allowing
$\Gamma$ to be a finite-index subgroup of $\Sp(2g,\ZZ)$ does not
change the argument. To fix notation, we choose a free $\ZZ$-module
$\WW$ of rank $2g$ equipped with a $\ZZ$-basis $e_1,\ldots,e_{2g}$,
and fix a standard skew-symmetric form by requiring
$\latt{e_i,e_{i+g}}=-1$ for $1\le i\le g$ and $\latt{e_i,e_j}=0$ if
$|i-j|\neq g$, so that the matrix of the form with respect to the
basis $\{e_i\}$ is
$\begin{pmatrix}0&-\Eins_g\\ \Eins_g&0 \end{pmatrix}$
(we denote the
$r\times r$ identity matrix by $\Eins_r$ throughout). Then
\[
\Sp(2g,\ZZ)=\left\{\gamma \in \GL(2g,\ZZ) \mid
   {}^t\gamma\begin{pmatrix}0&-\Eins_g\\ \Eins_g&0 \end{pmatrix}\gamma
   = \begin{pmatrix}0&-\Eins_g\\ \Eins_g&0 \end{pmatrix}\right\}.
\]

Choose a cusp $F_{g'}$ of rank $g'$ of $\HH_g$, for some $0\le g'<g$,
and put $g''=g-g'$. Recall that~$\HH_g$ is the cone of symmetric
$g\times g$ complex matrices $M$ with positive definite imaginary
part: cusps of rank $g'$ arise by allowing $\Imag M$ to be semidefinite
but requiring the radical $\rad M$ to be defined over $\QQ$ and of
dimension $g''$. Such a cusp thus corresponds to a rank~$g''$
isotropic sublattice $\XX$ of~$\WW$ up to the action of $\Gamma$, but
$\Gamma=\Sp(2g,\ZZ)$ acts transitively on such lattices. Therefore,
without loss of generality, we may take $F_{g'}$ to have stabiliser
\[
P(g')=\left\{\begin{pmatrix} A & 0 & B & m' \\ m & u & n & M \\ C & 0
& D & n'\\ 0 & 0 & 0 & {}^tu^{-1}
\end{pmatrix}\mid \begin{pmatrix} A&B \\ C&D
\end{pmatrix}\in \Sp(2g',\RR), \ u\in \GL(g'',\RR)\right\}\\
\]
in $\Sp(2g,\RR)$, by choosing $\XX=\ZZ e_{g'+1}+\ldots+\ZZ e_g\imic
\ZZ^{g''}$.

Next we consider the integral affine symplectic group
$\Tilde\Gamma_g^n$.
It is given by
\[
\Tilde\Gamma^n=\ZZ^{2gn}\ltimes \Sp(2g,\ZZ)
< \Tilde\Gamma_\RR^n=\RR^{2gn}\ltimes \Sp(2g,\RR)
<\GL(n+2g,\RR),
\]
and consists of elements $\Tilde\gamma$ of the form
\begin{equation}\label{eq:gammatilde}
\Tilde\gamma = \begin{pmatrix}1&a&b\\ 0&A_0&B_0\\ 0&C_0&D_0
\end{pmatrix}, \qquad \gamma = \begin{pmatrix} A_0&B_0 \\ C_0&D_0
\end{pmatrix}\in \Sp(2g,\RR), \ a,\,b\in M_{n\times g}(\RR)
\end{equation}
(cf.~\cite[Paragraph~(2.7)]{nam}).

The integral affine symplectic group acts on $\CC^{gn}\times\HH_g$, and
the Kuga variety, cf.~\cite[Equation~(3.4.1)]{nam}, is the quotient
\[
X_g^n:=
\Tilde\Gamma^n\backslash(\CC^{gn}\times\HH_g).
\]
It is the coarse moduli space of principally polarised abelian
varieties together with a point on the Kummer variety of the $n$-fold
cartesian self-product.

The stabiliser of $F_{g'}$ in $\Tilde\Gamma_{\RR}^n$ is
(cf.~\cite[Example~(2.8)]{nam})
\begin{equation}\label{eq:gammadash}
\Tilde{P}(g')=\left\{\begin{pmatrix}
\Eins_n & a' & a'' & b' & b''\\ 0 & A & 0 & B & m'\\ 0 & m & u & n & M\\ 0 &
C & 0 & D & n'\\ 0 & 0 & 0 & 0 & {}^tu^{-1}
\end{pmatrix}
\,\middle|\, 
\begin{array}{c}
\gamma' = \begin{pmatrix} A&B \\ C&D
\end{pmatrix}\in \Sp(2g',\RR),\\ \ u\in \GL(g'',\RR),\\ a',\,b'\in
M_{n\times g'}(\RR),\\ \ a'',\,b''\in
M_{n\times g''}(\RR)
\end{array}
\right\}
\end{equation}
where as before $M$, $m$ and $n$, and $m'$ and $n'$ are subject to the
symplecticity conditions, and its unipotent radical has centre
\[
\Tilde{U}(g')=\{u(b'',M)\mid M={}^tM\},\text{ where
}u(b'',M)=\begin{pmatrix}1&0&0&0&b''\\0&1&0&0&0\\0&0&1&0&M\\0&0&0&1&0\\0&0&0&0&1
\end{pmatrix}.
\]
Intersecting with $\Tilde\Gamma^n$ we obtain the group
\[
\Tilde\Upsilon^n=\Tilde U(g')\cap \Tilde P(g') = \{u(b'',M)\mid b''\in
\ZZ^{ng''},\, M={}^tM\in M_{g''\times g''}(\ZZ)\},
\]
which is identified with $\Sym^2(\XX^\vee)\times (\XX^\vee)^n$. (This
means a symmetric bilinear function on $\XX$ and $n$ linear functions
on $\XX$. In other words, $\Sym^1(\XX^\vee)$ is $\XX^\vee$. It is what
is called $B(X)\oplus (X^*)^n$ in \cite{fc}.)
  
To obtain the partial compactification at the cusp $F_{g'}$, we first
take the partial quotient by~$\Tilde\Upsilon^n$. For this we use the
Siegel domain realisation of $\HH_g$: for $\tau\in\HH_g$ we write
\[
\tau=\begin{pmatrix}\tau'&\omega\\ {}^t\omega&\tau''
\end{pmatrix}
\]
with $\tau'\in \HH_{g'}$, $\omega\in M_{g'\times g''}(\CC)$ and
$\tau''\in M_{g''\times g''}^{\sym}(\CC)$, and then
\[
\HH_g\imic\cD_{g'}:=\{(\tau',\omega,\tau'')
  \mid\Imag \tau'' -(\Imag {}^t\omega)(\Imag \tau')^{-1}(\Imag \omega) >0\}.
\]
Then $M\in \Sym^2(\XX)$ acts by translations in the imaginary
directions in $M_{g''\times g''}^{\sym}(\CC)$, so near this boundary
$X_g^n$ is covered by
\begin{equation}\label{eq:siegeldomain}
  \cD_{g'} \times \CC^{ng'}\times (\CC^*)^{ng''}\subset \HH_{g'}\times
  \CC^{g'g''}\times (\CC^*)^{g''\times g''}_{\sym} \times
  (\CC^{g'}\times (\CC^*)^{g''})^n
\end{equation}
where the $((\CC^*)^{g''})^n$ term is $\CC^{ng''}/(\XX^\vee)^n$, given
by $b''$ acting by translations in the imaginary directions in
$\CC^{ng''}=\XX^\vee\otimes\CC$.

Now we compactify by replacing the torus part $(\CC^*)^{g''\times
  g''}_{\sym}\times ((\CC^*)^{g''})^n$ (that is, the quotient
$\CC^{({g''}^2+(2n+1)g'')/2}/\Tilde\Upsilon^n$) by a suitable torus
embedding $\Temb(\Sigma(g'))$ corresponding to a fan~$\Sigma(g')$ in
$(\Sym^2(\XX^\vee)\times (\XX^\vee)^n)\otimes \RR$. We first define
two (non-polyhedral) cones, in $\Sym^2(\XX^\vee)$ and~$\Sym^2(\XX^\vee)\times (\XX^\vee)^n$ respectively,
following~\cite{fc}.

The first cone is $C(\XX)\subset \Sym^2(\XX^\vee\otimes\RR)$, which is
defined to be the cone of positive semi-definite symmetric bilinear
forms $b$ on $\XX^\vee\otimes\RR$ with rational radical $\rad
b$. Equivalently, $C(\XX)$ is the cone generated by the rank-$1$ forms
in $\Sym^2(\XX^\vee\otimes\QQ)$. To construct a toroidal
compactification of $\cA_g$ one must give decompositions of these
cones for each $g_i<g$ into fans invariant under the stabiliser in
$\Gamma$ of the cusp $F(g')$. There are many ways to do that, but
using reduction theory of quadratic forms (which guarantees the
compatibility) supplies two important ones, the first and second
Voronoi decomposition.

The second cone is $\Tilde C(\XX)\subset (\Sym^2(\XX^\vee)\oplus
(\XX^\vee)^n)\otimes\RR$, given by
\[
\Tilde C(\XX)=\{(b,\ell_1,\ldots,\ell_n)\in
\Sym^2(\XX^\vee\tensor\RR)\oplus (\XX^\vee\tensor\RR)^n \mid
\ell_i|_{\rad b}=0\text{ for all }i\}.
\]
as described in \cite[Definition~VI.1.3]{fc}. The number $n$ does not
change so we suppress it in the notation, but observe that $\Tilde
C(\XX)$ implicitly depends on $n$ whereas $C(\XX)$ does not.

We then further decompose $\Tilde C(\XX)$ as in
\cite[Definition~VI.1.3]{fc}, for each $g'$, obtaining a collection of
fans $\Tilde\Sigma=\{\Tilde\Sigma(g')\mid 0\le g' <g\}$.  Provided
that we choose the fans compatibly for different cusps we obtain, by
the standard toroidal compactification procedure from~\cite{AMRT}, a
compactification $\Bar {X_g^n}^{\Tilde\Sigma}$ that is analytically
locally isomorphic to the product of a smooth space with a quotient of
$\Temb(\Tilde\Sigma(g'))$ by $\Tilde\Gamma^n\cap P(g')$.

Ultimately we shall choose $\Tilde\Sigma(g')$ to be regular, so that
$\Temb(\Tilde\Sigma(g'))$ is smooth.  Singularities will then arise when we
move beyond $\Tilde\Upsilon^n$ and take the quotient by the rest of
$\Tilde\Gamma^n$, which in our situation may have fixed points.

We need more, however, because we seek a Namikawa
compactification. These are alluded to in~\cite{fc}, but the
equidimensional condition~\cite[Definition~VI.1.3(v)]{fc}, which is not
used and therefore not examined in detail there, is crucial for us. To
be precise, we require a weaker version of equidimensional, which we
call equidimensional in codimension~1. For this, it
is enough if every ray in $\Tilde\Sigma(g')$ maps onto a cone of the
fans $\Sigma$ that define $\Bar\cA_g$, whereas equidimensionality
requires this for every cone in $\Tilde\Sigma(g')$ of any dimension
for every ray $\tau$.

These two conditions, smoothness and equidimensionality, are in general
opposed to one another. Choosing $\Tilde\Sigma$ to give smooth
covering spaces $\Temb(\Tilde\Sigma(g'))$ typically involves blowing
up, and thus instantly violates condition~(iii) of the Namikawa
compactification.

Therefore, to construct an appropriate $\Bar{X_g^n}$, we need a
slightly more indirect approach. Instead of taking a regular
decomposition straight away, we first choose a decomposition
$\Tilde\Sigma^\flat$ such that $\Temb(\Tilde\Sigma^\flat(g'))$ itself
has canonical singularities. We do this by extending the perfect cone,
or first Voronoi, compactification of $\cA_g$, which has this property
by construction.

\begin{proposition}\label{prop:perfectcone}
There exist a $\Gamma$-admissible collection $\Sigma^\flat$ of fans
$\Sigma^\flat(g')$, for $0\le g'\le g$, and a~$\Tilde\Gamma$-admissible collection $\Tilde\Sigma^\flat$ of fans
$\Tilde\Sigma(g')$ such that
  \begin{enumerate}[(i)]
    \item $|\Tilde\Sigma^\flat(g')|=\Tilde C(\XX)$ and
      $|\Sigma^\flat(g')|=C(\XX)$;
    \item $\Tilde\Sigma^\flat(g')$ is $\GL(\XX)\ltimes \XX^n$-admissible relative
      to $\Sigma^\flat(g')$, for each $g'$;
    \item $\Temb(\Tilde\Sigma^\flat(g'))$ has canonical singularities
      and $\Tilde\Sigma^\flat(g')$ is equidimensional in codimension~1
      over $\Sigma^\flat(g')$;
    \item $\sigma\times\{0\}\in \Tilde\Sigma^\flat(g')$ for every $\sigma\in\Sigma^\flat(g')$.
  \end{enumerate}
\end{proposition}

\noindent{\it Remark.\/} The uses of the word ``admissible'' in the
preamble to Proposition~\ref{prop:perfectcone} and in condition~(ii)
there are different. $\Gamma$-admissible refers to the property of a
collections of fans, one for each cusp, being compatible under
restriction (see~\cite{AMRT}), whereas $\GL(\XX)\ltimes
\XX^n$-admissible is a property of compatibility at each cusp
separately with the the projection map $X_g^n\to \cA_g$, defined
in~\cite[Definition~VI.1.3]{fc}.

\noindent{\it Proof.\/} We take $\Sigma^\flat(g')$ to be defined by
the perfect cone decomposition of $C(\XX)$. This is the same as taking
the the cones of $\Sigma^\flat(g')$ to be the cones on the faces of
the convex hull of the rank-$1$ forms in the closure of $C$ with
rational radical, by~\cite{bc}. It is known to give an admissible
decomposition and a polyhedral fundamental domain for the action of
$\GL(\XX)$: see, for example, \cite[Ch.~8]{nambook} and the references
there.

We can extend this to a decomposition of $\Tilde C(\XX)$ by taking the
convex hull of all $(b,\ell)\in\Tilde C(\XX) \cap \Sym^2(\XX)\oplus
(\XX^\vee)^n$ with $\rank b=1$.  From the description of the action of
$\GL(\XX)\ltimes \XX^n$, given for example in \cite[VI.1.1]{fc}, it
follows immediately that $\GL(\XX)\ltimes \XX^n$ acts on $\Tilde
C(\XX)$ with a polyhedral fundamental domain.

If $q=(b;\ell_1,\ldots,\ell_n)=(b;\ell_j)\in \Tilde C(\XX)$ then $b\in
C(\XX)$ so we can write $b=\sum \lambda_i r_i$, where ${r_i=\xi_i{}^t\xi_i}$ is of rank~$1$, $\xi_i\in \XX$, and $\lambda_i\in
\RR_+$. Since $\ell_j|_{\rad b}$ we may write
$\ell_j=b(\mu_j,\bullet)$ with $\mu_j\in \rad b$ (notice that $\rad b
\supseteq \bigcup_i \rad r_i$, with equality if the $r_i$ are linearly
independent). Hence
\[
q=(b;\ell_j)=(b;b(\mu_j,\bullet))
=(b;\sum_i\lambda_ir_i(\mu_j,\bullet))=\sum_i\lambda_i(r_i;r_i(\mu_j,\bullet)),
\]
so that the cone over the convex hull of the integral $q$ with
rank~$1$ quadratic part $b$ is indeed~$\Tilde C(\XX)$, and every
integral point of $\Tilde C(\XX)$ is in the convex hull.

The first of these conditions shows that the definition of
$\Tilde\Sigma^\flat$ does give an admissible collection of fans. That
is, it is $\Tilde\Gamma^n$-invariant, and chosen for different cusps
$F_{g'}$ so as to be compatible with restriction (Siegel
$\Phi$-operator) to adjacent cusps. This therefore yields a
compactification $\Bar{X_g^n}^\flat$.

The second condition shows that the covering spaces
$\Temb(\Tilde\Sigma^\flat(g'))$ have canonical singularities.

Finally, $\Bar{X_g^n}^\flat$ is a Namikawa compactification because
there are no rays of $\Tilde\Sigma^\flat(g')$ in the interior of
$\Tilde C(\XX)$, and the ray spanned by $q=(b;\ell_j)$ projects onto
the ray spanned by $b$. \qed\smallskip
\medskip

Shepherd-Barron showed in~\cite{s-b} (see also the
correction~\cite{as-b}) that the perfect cone compactification of
$\cA_g$ has canonical singularities for $g \ge 5$. We are not
constrained to use a specific compactification: rather, we choose a
suitable one, as in~\cite[Section~5]{Tai}. We choose smooth
subdivisions $\Tilde\Sigma^\sharp(g')$ of the fans $\Tilde\Sigma^\flat(g')$
(that is, toric resolutions of $\Temb(\Tilde\Sigma^\flat(g'))$) in a~$\Tilde\Gamma$-equivariant way, and denote the resulting
compactification by $\Bar{X^n_g}^\sharp$. This is of course no longer a
Namikawa compactification, nor is it smooth in general since
$\Tilde\Gamma$ is not neat. However, we have the following easy
consequence of Lemma~\ref{lem:covercanonical}.

\begin{corollary}\label{cor:coverenough}
Suppose that $\Bar{X^n_g}^\sharp$ has canonical singularities and that
the action of $P(g')\cap\Tilde\Gamma$ on
$\Temb(\Tilde\Sigma^\sharp(g'))$ has no quasireflections. Then
$\Bar{X_g^n}^\flat$ has canonical singularities.
\end{corollary}

\noindent{\it Proof.\/} It is enough to apply
Lemma~\ref{lem:covercanonical} to the resolutions
$\Temb(\Tilde\Sigma^\sharp(g'))\to
\Temb(\Tilde\Sigma^\flat(g'))$ for each $g'$. \qed\smallskip

\section{Quotient singularities}\label{sect:quotient}

In this section we shall verify the conditions of Corollary~\ref{cor:coverenough}
for the values of $g$ and $n$ that concern us.

Recall that any toroidal compactification $\Bar\cA_g$ of $\cA_g$ comes
with a surjective map to the Satake compactification $\cA_g^*$, and
that the Satake compactification has a stratification
\[
\cA_g^*=\cA_g\sqcup\cA_{g-1}\sqcup\ldots\sqcup\cA_1\sqcup\cA_0
\]
where $\cA_0$ is a point. Hence any compactification of
$X_g^n$ that dominates a Namikawa compactification has a map to
$\cA_g^*$ and in particular there is a natural map $\pi\colon
\Bar{X_g^n}^\sharp\to \cA_g^*$.

If $p\in\Bar{X_g^n}^\sharp$ is such that $\pi(p)\in\cA_{g'}\subset
\cA_g^*$, then near $p$ the toroidal compactification is a quotient of
an invariant open subset
\begin{equation}\label{eq:productcover}
  \cD^\sharp\subset \HH_{g'}\times  \CC^{g'g''}\times \CC^{g'n}\times \Temb(\Tilde\Sigma^\sharp)
\end{equation}
(where $g'+g''=g$) by an action of $\Tilde P(g')\cap \Tilde\Gamma^n$
that preserves the product structure, extending the
decomposition~\eqref{eq:siegeldomain}. To determine the singularity at
$p\in\Bar{X_g^n}^\sharp$ we therefore have to examine the action of
the stabiliser $G$ of a preimage $\tilde p\in \cD^\sharp$ of $p$ on
the tangent space $T_{\cD^\sharp,\tilde p}$.

We recall some basic facts from \cite{Tai} and \cite{YPG}. Suppose
that $\rho\colon G \to \GL(m,\CC)$ is a finite-dimensional
representation of a finite group $G$, and suppose that $h\in G$
and that $\rho(h)$ has order $k>1$. If the eigenvalues of
$\rho(h)$ are $\zeta^{a_1},\ldots,\zeta^{a_m}$ (where
$\zeta=e^{2\pi i/k}$ is a primitive $k$-th root of unity and $0\le
a_i<k$) then the \emph{Reid-Tai sum} of $h$, also called the
\emph{age} of $h$, is
\[
\RT(h)=\sum_{i=1}^m \frac{a_i}{k}.
\]
The RST (Reid--Shepherd-Barron--Tai) criterion states that if
$\rho(G)$ has no quasireflections then the quotient $\CC^m/\rho(G)$
has canonical singularities if and only if $\RT(h)\ge 1$ for every
$h\in G$.

If $\rho(h)$ is a quasireflection then exactly one of the $a_i$ is
non-zero, so $\RT(h)<1$. It follows that in any case if $\RT(h)\ge 1$
for every $h\in G$ then the quotient has canonical singularities.

We apply this to the action of $\tilde\gamma\in \Tilde P(g')$. We use
the block decomposition given in~\eqref{eq:gammatilde}
and~\eqref{eq:gammadash} and the notation for submatrices in the rest
of this section is taken from there.

To check that the singularity at $p$ is canonical it is enough to
verify that $G$ contains no quasireflection on the tangent space
(which we need to do anyway in order to apply
Lemma~\ref{lem:covercanonical}) and that $\RT(\Tilde\gamma)\ge 1$ for
any nontrivial $\Tilde\gamma\in G$.

Note that the decomposition~\eqref{eq:productcover} is $G$-invariant,
so that $\RT(\Tilde\gamma)$ is the sum of the age of $\Tilde\gamma$
restricted to each factor.

\begin{proposition}\label{prop:zerobound}
If $g''=0$ then $\Bar{X_g^n}$ has a canonical singularity at $p$ and
the stabiliser $G$ of $\tilde p$ has no quasireflections, unless $g=2$
and $n\le 2$, or $g=3$ and $n=1$.
\end{proposition}

\noindent{\it Proof.} Recall that we are assuming $g\ge 2$ anyway. If
$g''=0$ then the local cover in~\eqref{eq:productcover} becomes~${\cD^\sharp=\HH_g\times \CC^{gn}}$, which is just the covering space of
$X^n_g$. But $X^n_g$ has canonical singularities by
\cite[Proposition~10.3]{ma}, and there are no quasireflections by
\cite[Lemma~7.1]{ma}, for these values of $g$ and $n$. \qed\smallskip

In view of Proposition~\ref{prop:zerobound}, we may assume for the
rest of this section that $g''>0$.

In order to use the results of~\cite{Tai} we need to verify the
condition that $\tilde\gamma$ should act trivially on each cone
of the fan (see \cite[p.~438]{Tai}). 

\begin{lemma}\label{lem:trivialsigmaaction}
If $\tilde\gamma$ fixes a cusp then it acts trivially on each cone of
$\Tilde\Sigma^\sharp(g')$.
\end{lemma}

\noindent {\it Proof.}  The eigenvectors in $\Sym^2(V)$ are the
rank~$1$ forms $f_if_j$, where $f_i\in\XX^\vee\tensor\RR$ are the real
eigenvectors of $u$ (if any) and thus there are no eigenvectors in the
interior of the cone $\Tilde C$.  Now the result follows: for if
$\tilde\gamma$ acts nontrivially on $\sigma\in\Tilde\Sigma^\sharp(g')$
then it acts on some subset $\{\rho_1,\ldots,\rho_k\}$ of the rays
spanning $\sigma$ by a free permutation, and then by the Brouwer
fixed-point theorem it has an eigenvector in the closed cone spanned
by $\{\rho_1,\ldots,\rho_k\}$; but this eigenvector is not a generator
of any of the $\rho_i$, so it is in the interior of $\Tilde
C$. \qed\smallskip

\begin{proposition}\label{prop:unontrivial}
Suppose that $p\in \Bar{X^n_g}$ and $\pi(p)\in \cA_{g'}$, and that
$\Tilde\gamma\neq\Eins_{n+2g}$ belongs to the stabiliser $G$ of
$\tilde p\in\cD^\sharp$ in $\Tilde P(g')\cap \Tilde\Gamma^n$. Suppose
that $u\neq \pm\Eins_{g''}$, or that $u=-\Eins_{g''}$ and
$(g'',n)\neq(1,1)$. Then the action of $\Tilde\gamma$ on
$T_{\cD^\sharp,\tilde p}$ has $\RT(\Tilde\gamma)\ge 1$, and in
particular it is not a quasireflection.
\end{proposition}

\noindent {\it Proof.}  For this, it is enough to look at the
$\Temb(\Tilde\Sigma^\sharp)$ factor. This is a toric variety with
torus $\TT$ whose lattice of 1-parameter subgroups (the dual of the
character lattice) is $\Sym^2(\XX^\vee)\times(\XX^\vee)^n$.

If the eigenvalues of $\tilde\gamma$ on $V=\XX^\vee\otimes \CC$, which
are the eigenvalues of $u$, are $\mu_1,\ldots,\mu_{g''}$ then the
eigenvalues on $\Sym^2(V)\times V^n$ are $\mu_i\mu_j$ and $n$ copies of
each $\mu_i$. The eigenvalues of~$\tilde\gamma$ on the toric boundary
component containing $\tilde p$ belongs include $n$ copies of each
$\mu_i$ and, as in \cite[Lemma~5.2]{Tai}, all the $\mu_i\mu_j$ that
are different from~$1$. Since $u\in\GL(g'',\ZZ)$ is of finite order, its
eigenvalues include all the primitive $d$-th roots of unity for some
degree, and that gives~$\RT(\tilde\gamma|_V) \ge 1$ on the $V$ unless
$d=1$ or $d=2$.

If $d=1$ then $u=\Eins_{g''}$. If $d=2$ then $\mu_i=\pm 1$ and we may
assume $\mu_1=-1$: then if $g''>1$ we either have $\mu_2=-1$ so again
$\RT(\tilde\gamma|_V) \ge 1$, or $\mu_2=1$ and then the eigenvalue
$\mu_1\mu_2=-1$ occurs on $\Sym^2(V)$. If $g''=1$ then the eigenvalue
$\mu_1=-1$ occurs $n$ times on $V^n$, so $\RT(\tilde\gamma|_{V^n}) \ge
1$ unless $g''=n=1$.  \qed\smallskip
\medskip

Next, we examine the action of $\tilde\gamma$ on the $\CC^{ng'}$
factor. Because of Proposition~\ref{prop:unontrivial} we may assume
that $u=\epsilon\Eins_{g''}$ with $\epsilon=\pm 1$.

For any $r\leq g$, we let $M_{r\times g}(\CC)^*$ be the set of
matrices of rank $r$ in $M_{r\times g}(\CC)$. Then the Grassmannian
$\Gr(r, g)$, of $r$-dimensional linear subspaces in $\CC^g$, is
$M_{r\times g}(\CC)^*/\GL(r,\CC)$, with~$\GL(r,\CC)$ acting by right
multiplication. Since $\CC^{ng}$ is identified with $M_{n\times
  g}(\CC)$ by the choice of basis~$e_1,\ldots,e_{2g}$, we may regard
$\HH_g\times \CC^{ng}$ as a subset of subset of $\Gr(g, n+2g)$ by
sending an element $(\tau,Z)\in \HH_g\times \CC^{ng}$ to the
equivalence class of block matrices:
\[
\begin{bmatrix}
Z\\
\tau\\
\Eins_g
\end{bmatrix} \in M_{(n+2g) \times g}(\CC)/\GL(g,\CC).
\]
Recall that in this representation a boundary component of $\HH_g$ is
a subset of the closure $\Bar\HH_g$ of~$\HH_g$ in $\Gr(g,2g)$, so this
description extends to the boundary. So the image of $\Tilde p$ in
$\Bar\HH_g\times\CC^{ng}$ may also be written in this way, with
$\tau\in F_{g'}\subset \Bar\HH_g$ given by
\[\tau = 
\begin{pmatrix}
\tau' & \omega\\
{}^t\omega & \tau''
\end{pmatrix}
\]
where $\tau' \in \HH_{g'}$, $\omega \in M_{g'\times g''}$ and $\tau'' \in
\Sym^2(\XX^\vee\otimes \CC)=M^{\sym}_{g''\times g''}(\CC)$.

Then the action of $\tilde\gamma$ is given (notation
from~\eqref{eq:gammatilde}) by
\[\tilde\gamma \cdot
\begin{bmatrix}
Z \\ \tau \\ \Eins_g
\end{bmatrix} = 
\begin{bmatrix}
Z+a\tau+b\\
\begin{matrix}
A\tau'+B & Aw+m'\\
m\tau' + \epsilon{}^tw + n & mw + \epsilon \tau''+M\\
C\tau' + D & Cw + n'\\
0 & \Eins_{g''}
\end{matrix} 
\end{bmatrix}.
\]
Because $\tau$ is preserved by $\gamma$, this simplifies to
\[\tilde\gamma \cdot
\begin{bmatrix}
Z \\ \tau \\ \Eins_g
\end{bmatrix} = 
\begin{bmatrix}
(Z+a\tau+b) \cdot N\\
\tau\\
\Eins_g
\end{bmatrix}
\]
where 
\[
N = 
\begin{pmatrix}
(C\tau' + D)^{-1} & -(Cw + n')(C\tau' + D)^{-1}\\
0 & \Eins_{g''}
\end{pmatrix},
\]
and therefore the action of $\tilde\gamma$ on the tangent space to
$\CC^{ng}$ at $Z$ is by right multiplication by $N$.

\begin{lemma}\label{lem:evalsonnpart}
The eigenvalues of $\tilde\gamma$ on $\CC^{ng}$ are exactly the
eigenvalues of $N$, which are $1$ and the eigenvalues of $(C\tau' +
D)^{-1}$. Moreover, $\gamma'$ fixes $\tau'$, i.e.\ $(A\tau' +
B)(C\tau' + D)^{-1} = \tau'$.
\end{lemma}

\noindent{\it Proof.} Immediate from the discussion
above. \qed\smallskip
\medskip

The next lemma and its corollaries apply to the usual action of
$\Sp(2g,\ZZ)$ on $\HH_g$, for arbitrary~$g$. If
$\gamma=\begin{pmatrix}A_\gamma & B_\gamma\\ C_\gamma & D_\gamma
\end{pmatrix}
\in \Sp(2g,\ZZ)$ then denote the eigenvalues of $\gamma$ by
$\Lambda=(\lambda_1\dots,\lambda_g\}$ and
${\Bar\Lambda=(\bar\lambda_1,\dots,\bar\lambda_g)}$: we will regard such
a sequence as a diagonal matrix.

\begin{lemma}\label{lem:evalsoccur} If $\gamma$ fixes $\tau\in\HH_g$, then
  $\gamma$ is diagonalisable. Moreover $C_\gamma\tau + D_\gamma$ is
  diagonalisable and has $\bar\Lambda$ as sequence of eigenvalues.
\end{lemma}

\noindent{\it Proof.}  For every $\beta\in \Sp(2g,\RR)$ and $\tau\in
\HH_g$ we set $J(\beta,\tau) = C_\beta\tau+D_\beta$. Then $J$ is a
cocycle, i.e.\ $J(\beta_1\beta_2,\tau)=J(\beta_1,
\beta_2\tau)J(\beta_2, \tau)$ for every
$\beta_1,\,\beta_2\in\Sp(2g,\RR)$ and $\tau\in\HH_g$.
 
It is well known,  cf.~\cite[Lemma~4.1]{Tai}, that there exists
$\alpha \in \Sp(2g,\RR)$ such that
\[
\alpha\tau=i \Eins_g\qquad\text{and}\qquad
 \alpha\gamma\alpha^{-1}= \begin{pmatrix} \delta_1&\delta_2
   \\ -\delta_2&\delta_1 \end{pmatrix}
\]
with $\delta_1$, $\delta_2$ real diagonal matrices and
$\delta_1+i\delta_2 \in \Ut(g, \CC)$.

Obviously $\gamma$ is diagonalisable, with eigenvalues $\Lambda =
\delta_1+i\delta_2$ and $\Bar\Lambda=\delta_1-i\delta_2$.  Now using the
cocycle property we have
\begin{align*}
  \Bar\Lambda &=J(\alpha\gamma\alpha^{-1}, i\Eins_g)\\
  &=J(\alpha, \gamma\alpha^{-1},i\Eins_g)J(\gamma\alpha^{-1}, i\Eins_g)\\
  &= J(\alpha, \gamma\cdot\tau)J(\gamma, \tau) J(\alpha^{-1},
  i\Eins_g)
\end{align*}
Now, since $\gamma$ fixes $\tau$ and $J(\alpha^{-1}, i\Eins_g)=
J(\alpha, \tau)^{-1}$, we get that $\Bar\Lambda$ and $J(\gamma,\tau)$
are conjugate.  \qed\smallskip
\medskip

We have the following immediate corollary.

\begin{corollary}\label{cor:nononeeval} If $\gamma \in \Sp(2g,\ZZ)$ is
  nontrivial and fixes $\tau\in\HH_g$, then $C_\gamma\tau + D_\gamma$ has an
  eigenvalue that is not $1$.
\end{corollary}

Now we return to the singularities of $\Bar{X_g^n}$.
  
\begin{proposition}\label{prop:onebound}
If $g''\neq 0$ then $\Bar{X_g^n}$ has a canonical singularity at $p$,
as long as $g+n\ge 6$.
\end{proposition}

\noindent{\it Proof.} As before, we take $\Tilde\gamma\in \Tilde
P(g')$, fixing a point $\tilde p$, and write it using the block
decomposition given in~\eqref{eq:gammatilde} and~\eqref{eq:gammadash}.
Again, because of Proposition~\ref{prop:unontrivial} we may assume
that $u=\epsilon\Eins_{g''}$, where~${\epsilon=\pm 1}$ and $\epsilon=1$
unless $g''=n=1$.

For any $g$, if $\gamma'=\Eins_{2g'}$ and $u=\Eins_{g''}$ then
$\tilde\gamma\in \tilde U(g')$ (see~\cite[Example 2.8]{nam}) and acts
trivially at the boundary $F_{g'}$. In particular this holds if
$g'=0$, unless $g''=n=1$ but then $g=1$ which is excluded.

If $\gamma'=-\Eins_{2g'}$ and $\epsilon=1$, or $\gamma'=\Eins_{2g'}$ and
$\epsilon=-1$, then there are $g-1$ eigenvalues $\lambda_i\epsilon=-1$
on the $\CC^{g'g''}$ factor, giving $\RT(\tilde\gamma)>1$.

If $\gamma'=-\Eins_{2g'}$ and $\epsilon=-1$ then the eigenvalues on the
$\Temb(\Tilde\Sigma^\sharp)$ factor include $n$ copies of $\epsilon$,
and, by Corollary~\ref{cor:nononeeval}, there are also $n$ copies of
$-1$ occurring on the $\CC^{ng}$ factor. Even for $n=1$, this gives
$\RT(\gamma)\ge 1$.

Therefore we may assume that $\gamma'\neq \pm\Eins_{2g'}$, and thus
does not act trivially on $\HH_{g'}$.

If $g'\ge 5$ and $\gamma'\neq \pm \Eins_{2g'}$ then the contribution to
$\RT(\gamma)$ from $\gamma'$ acting on the $\HH_{g'}$ factor
(the $F_{g'}$ factor) is already at least~$1$, and $\cA_{g'}$ itself
has canonical singularities: this is~\cite[Lemma~4.5]{Tai}.

If $g'<5$ then the eigenvalues not coming from the action of $\gamma'$
on $\HH_{g'}$ include $g''$ copies of~$\epsilon\lambda_i$ on the
$\CC^{g'g''}$ factor and, by Corollary~\ref{cor:nononeeval}, a further
$n$ copies of $\epsilon\lambda_i^{\pm1}$ on the $\CC^{ng}$ factor. If
the order of $\tilde\gamma$ on the tangent space at $\tilde p$ is $d$,
then some $\lambda_i$ is a nontrivial $d$-th root of unity and so this
gives a contribution of at least $(n+g'')\frac{1}{d}$ to $\RT(\tilde\gamma)$.

Moreover, according to~\cite[Lemma~4.4]{Tai} we may assume $d\le 6$,
and in each case the action of $\gamma'$ contributes at least
$\frac{g'}{d}$ to $\RT(\tilde\gamma)$, so in any case we have
$\RT(\tilde\gamma)\ge (g'+n+g'')\frac1d \ge \frac{n+g}6$. So if
$n+g\ge 6$ we are done.~\qed\smallskip
\medskip

The condition $n+g\ge 6$ cannot be strengthened: it is needed if
$g'=1$ and $d=6$.

Theorem~\ref{thm:goodcompact} now follows immediately.

\section{Slope of $\cA_g$}\label{sect:slope}
We assume throughout that $g>1$, since the case of $g=1$ can be
reduced to the case of $\cM_{1,n+1}$, which is solved in~\cite{belo}
and~\cite{bf}.

We shall construct differential forms by using Siegel modular forms, so
we begin with some elementary definitions concerning them.

\begin{definition}\label{def:siegelmodform}
A \emph{modular form of weight $k$} is a
holomorphic function $f \colon\HH_g\to \CC$ on the Siegel upper
half-plane
\[
\HH_g=\{ Z\in M_{g\times g}(\CC)\mid Z={}^tZ,\ \Imag Z >0\}
\]
such that
\[
f(\gamma\cdot \tau)= \det (C\tau+D)^k f(\tau) \text{ for any } \gamma\in
\Sp(2g,\ZZ).
\]
\end{definition}
Note that we need no extra condition at infinity when $g>1$.

A Siegel modular form has a Fourier expansion
\[
f(\tau) =\sum_T a(T) \exp(\pi i \tr(T\tau))
\]
where the sum runs over all even integral symmetric matrices $T$.

\begin{definition}\label{def:vanishing}
If $f$ is a Siegel modular form, the \emph{vanishing of $f$ at the boundary} is
\[
b:=\half \min\{x^tTx\mid a(T)\neq 0,\ x\in \ZZ^g\smallsetminus\{0\}\}.
\]
If $b>0$, i.e.\ if $a(T)\neq 0$ implies $T>0$, we say that $f$ is a \emph{cusp form}.
\end{definition}

We recall, mainly from~\cite{mumfordimag}, some facts about
$\cA_g=\Sp(2g,\ZZ)\backslash \HH_g$ and its compactifications.

Modular forms of weight~$1$ determine a $\QQ$-line bundle $L$, the
Hodge line bundle. The Satake compactification $\cA_g^*$ is Proj of
the ring of modular forms, and the Mumford partial compactification
$\cA'_g$ is the blow-up of $\cA_g\sqcup\cA_{g-1}\subset\cA_g^*$ along $\cA_{g-1}$. Every toroidal compactification of~$\cA_g$ dominates
$\cA_g^*$ and contains $\cA'_g$ as a Zariski open subset: the toroidal
compactifications differ from one another only above the deeper strata
$\cA_{g'}$ for $g'<g-1$.

If $g\geq 2$ then $\Pic(\cA_g')\otimes\QQ=\QQ\lambda\oplus \QQ\delta$,
where $\lambda$ is the class of $L$ and $\delta$ is the class of the
boundary divisor $\Delta_{\cA'}$, the proper transform of
$\cA_{g-1}\subset \cA^*_g$. This is proved
in~\cite[Corollary~1.6]{mumfordimag} for $g\ge 4$, and for $g=2$ and
$g=3$ it follows from the studies of the Chow ring in \cite{mumenum}
and \cite{faber} respectively: see~\cite{Hu}.

\begin{definition}\label{def:slope}
The \emph{slope} of an effective divisor $E=a\lambda-b\delta$ on
$\cA'_g$ with $a,\,b> 0$ is defined to be $s(E)=a/b$. In particular,
if $E$ is the zero divisor of a cusp form $f$ of weight $k$ and
vanishing order $b$ then $s(E)=k/b$.
\end{definition}

The Kodaira dimension of the Kuga varieties $X_g^n$ is related to
the transcendence degree of the field generated by cusp forms of slope
less or equal to $g+n+1$: see for example~\cite[Theorem~1.3]{ma}.

\begin{definition}\label{def:minslope}
The \emph{minimal slope} $s_{\min}(g)$ is the infimum of the slopes of
all effective divisors on $\cA'_g$.
\end{definition}

An upper bound $s_{\min}(g)$ is provided in small genera by the
Andreotti-Mayer divisor $N_0$, the locus of principally polarised
abelian varieties with singular theta divisor~\cite{anma}. The divisor~$N_0$ has two components: $\Theta_{\Null}$, the locus where the theta
divisor has a singular point of order~$2$, and~$N'_0$, the locus where
the theta divisor has a singular point not of order~$2$. The classes
of~$\Theta_{\Null}$, for $g\ge 1$, and $N'_0$, for $g\ge 4$, are
computed in~\cite{mumfordimag}:
\begin{align*}
[\Theta_{\Null}] &= 2^{g-2}(2^g+1)\lambda - 2^{2g-5}\delta. \\
[N'_0]&=\left(\frac{(g+1)!}4+\frac{g!}{2}-2^{g-3}(2^g+1)\right)\lambda
-\left(\frac{(g+1)!}{24}-2^{2g-6}\right)\delta.
\end{align*}
For $g\le 3$ the minimal slope is achieved at $[\Theta_{\Null}]$,
giving the values $s_{\min}(1)=12$, $s_{\min}(2)=10$ and $s_{\min}(3)=9$.

For $g=4$ we have $s_{\min}(4)=8$, achieved by $s([N_0'])$, and for
$g=5$ we have $s_{\min}(5)=54/7$, also achieved by $s([N'_0])$: 
see \cite{sm} and \cite{fgsv}, respectively.

For $2\le g\le 4$, the divisors that minimise the slope are rigid.

For $g=6$ we have $s_{\min}(6)\le 7$: see~\cite{dssm}. However, the fact
that $s_{\min}(6)\le s([N'_0])=550/73<8$ will suffice for our purposes.

For $g\ge 7$ we have $s_{\min}(g)<g+1$ by~\cite{Tai} and \cite{mumfordimag}.

\begin{theorem}\label{thm:computekappa}
Suppose that $g\ge 2$ and $\Bar{X_g^n}$ is a Namikawa compactification of $X_g^n$
with canonical singularities. Then
\begin{enumerate}
\item $\kappa(X_g^n)=\half g(g+1)$ if $s_{\min}(g)<g+n+1$;
\item $\kappa(X_g^n)=0$ if $s_{\min}(g)=s(D)=g+n+1$ and $D$ is rigid;
\item $\kappa(X_g^n)=-\infty$ if $s_{\min}(g)>g+n+1$ (even if the
  singularities are not canonical).
\end{enumerate}
\end{theorem}

\noindent{\it Proof:\/} The first case (what one might call relatively
general type) follows easily from~\cite[Proposition~9.2]{ma}. Pulling
back along $f\colon \Bar{X_g^n}^\sharp\to \Bar\cA_g$, this
implies $kK_X\ge f^*(k(g+n+1)L-k\Delta_\cA)$ for sufficiently
divisible~$k$. So it is enough to show that the
$\QQ$-divisor $(g+n+1)L-\Delta_\cA$ is big: however, since it has
slope strictly greater than $s_{\min}(g)$ it is in the interior of the
effective cone and can therefore be written as the sum of an effective
divisor and an ample divisor.

In the second case, the same argument shows that $K_X$ is effective,
since it dominates the pullback of the effective divisor
$(g+n+1)L-\Delta_\cA$. Therefore $\kappa(X_g^n)\ge 0$. On the other
hand, if some multiple of $K_X$ moves, then so does some multiple of
$(g+n+1)L-\Delta_\cA$, which is to say that some multiple of $D$
moves, but $D$ is rigid.

In the third case, if $K_X\ge 0$ then $f_*(K_X)\ge 0$, but
$s(f_*(K_X))=g+n+1<s_{\min}(g)$. So $K_X$ is not effective, and
$\kappa(X_g^n)=-\infty$. \qed\smallskip

Theorem~\ref{thm:kdim} follows immediately from this.
\smallskip

The authors would like to thank G.~Farkas, S.~Ma and
N.I.~Shepherd-Barron for useful discussions on this topic, and the
referee for asking about the relation with $\cM_{g, n}$.

\bigskip
  \footnotesize

\noindent{}Flora Poon\\ \textsc{Mathematics Division, National Center
  for Theoretical Sciences,}\\  \textsc{No.~1, Sec.~4, Roosevelt~Rd., Taipei City, Taiwan.}\par\nopagebreak
\noindent\texttt{wkpoon@ncts.ntu.edu.tw}
\medskip

\noindent{}Riccardo Salvati Manni\\ \textsc{Dipartimento di Matematica ``Guido
    Castelnuovo'',}\\ \textsc{Universit\`a di Roma ``La
  Sapienza'',}\\  \textsc{P.le Aldo Moro 5, 00185 Roma, Italy}\par\nopagebreak
\noindent{}\texttt{salvati@mat.uniroma1.it}
\medskip

\noindent{}Gregory Sankaran\\ \textsc{Department of Mathematical
  Sciences,}\\ \textsc{University of Bath, Bath BA2 7AY,
  UK}\par\nopagebreak
\noindent{}\texttt{G.K.Sankaran@bath.ac.uk}

\end{document}